\documentclass[12pt]{article}
\usepackage{amsmath}
\usepackage{amssymb}
\usepackage[ctr1,eng]{myarticl}

\newcommand{\antip }{\kappa }
\newcommand{\Aut }{\mathrm{Aut}}

\newcommand{\cB }{\mathcal{B}}

\newcommand{\cG }{\mathcal{G}}

\newcommand{\cop }{^{\mathrm{cop}}}
\newcommand{\copr }{\varDelta }
\newcommand{\End }{\mathrm{End}}

\newcommand{\lact }{.}
\newcommand{\lacti }{.}
\newcommand{\lactZ }{\triangleright }

\newcommand{\lcoa }{\delta }
\newcommand{\lcoaZ }{\gamma }
\newcommand{\lcoai }{\delta }
\newcommand{\ldif }{y^{\scriptscriptstyle \mathrm{L}}}

\newcommand{\Ndbasis }{\boldsymbol{\mathrm{e}}}

\newcommand{\ndN }{\mathbb{N}}

\newcommand{\ndZ }{\mathbb{Z}}
\newcommand{\op }{^\mathrm{op}}
\newcommand{\opcop }{^{\mathrm{op},\mathrm{cop}}}
\newcommand{\ot }{\otimes }

\newcommand{\paar }[2]{\langle #1,#2\rangle }

\newcommand{\proots }{\boldsymbol{\Delta }^+}

\newcommand{\qnum }[2]{(#1)_{#2}}

\newcommand{\rdif }{y^{\scriptscriptstyle \mathrm{R}}}

\newcommand{\roots }{\boldsymbol{\Delta }}

\newcommand{\YD }{Yetter--Drinfel'd }
\newcommand{\YDcat }{\mathcal{YD}}

\title{The Weyl--Brandt groupoid of a\\
Nichols algebra of diagonal type
\thanks{Supported by the European Community under a Marie Curie
Intra-European Fellowship}}
\author{I.~Heckenberger}

\begin{document}

\maketitle

\begin{abstract}
The theory of Nichols algebras of diagonal type is known to be
closely related to that of semisimple Lie algebras. In this paper
the connection between both theories is made closer. For any
Nichols algebra of diagonal type invertible transformations are
introduced, which remind one of the action of the Weyl group on
the root system associated to a semisimple Lie algebra. They give
rise to the definition of a Brandt groupoid. As an application an
alternative proof of classification results of Rosso,
Andruskiewitsch, and Schneider is obtained without using any
technical assumptions on the braiding.


\end{abstract}

\section{Introduction}

A method of Andruskiewitsch and Schneider \cite{a-AndrSchn00} for
the classification of pointed Hopf algebras contains as an
essential step the determination of Nichols algebras with certain
growth conditions. However even in the simplest case, when the
Nichols algebra is of diagonal type, there exists up to now no
complete answer to the latter problem. General assertions which
relate Nichols algebras to semisimple Lie algebras were proved for
example by Rosso \cite{a-Rosso98} and Andruskiewitsch and
Schneider \cite{a-AndrSchn00}. In both cases the introduction of
additional technical assumptions was necessary.

Kharchenko \cite{a-Khar99} proved that any Nichols algebra $\cB
(V)$ of diagonal type has a (restricted)
Poincar\'e--Birkhoff--Witt basis consisting of iterated
skew-commutators of basis vectors of $V$. Further, $\cB (V)$ has a
natural $\ndZ ^n$-grading, where $n$ is the rank of $\cB (V)$.
This gives rise to the definition of the ``root system'' of $\cB
(V)$, see Section \ref{sec-fini}. It is natural to look for a
structure which plays here a similar role as the Weyl group for
ordinary root systems. In Section \ref{sec-trafo} transformations
of Nichols algebras of diagonal type are introduced. They
naturally give rise to a Brandt groupoid structure, see Section
\ref{sec-WBg}, associated to $\cB (V)$. It will be called the
\textit{Weyl--Brandt groupoid of} $\cB (V)$ and denoted by $W(V)$.
If the Nichols algebra is of \textit{Cartan type} then $W(V)$ is
isomorphic to $G\times B$ for a group $G$ and a set $B$. If the
corresponding Cartan matrix $C$ is symmetrizable then $G$ is
isomorphic to the Weyl group associated to $C$ and $B$ is the
orbit of the ordered standard basis of $\ndZ ^n$ under the action
of the Weyl group. In the last section of this paper $W(V)$ is
used to give a new proof of the classification results of Rosso
and Andruskiewitsch and Schneider without using technical
conditions on the braiding. The relative simplicity of the proof
stems from the fact that Kharchenko's results and the use of
$W(V)$ allow to determine the degrees and heights of the
(restricted) Poincar\'e--Birkhoff--Witt generators without knowing
the defining relations of $\cB (V)$ explicitly.

Throughout this paper $k$ denotes a field of characteristic zero
and tensor products $\ot $ are taken over this field. Given an
algebra $\cB $, let $\cB \op $ denote $\cB $ with the opposite
product. Similarly, if $C$ is a coalgebra then $C\cop $ denotes
$C$ with the opposite coproduct. For the coproduct and the
antipode of a Hopf algebra the symbols $\copr $ and $\antip $ are
used. The coproduct of elements $a$ of a coalgebra is written in
the Sweedler notation $\copr (a)=a_{(1)}\ot a_{(2)}$. The set of
natural numbers not including 0 is denoted by $\ndN $ and we write
$\ndN _0=\ndN \cup \{0\}$.

The author would like to thank Professor H.-J.~Schneider for
helpful discussions on Nichols algebras during the ``International
Conference on Quantum Groups'' dedicated to the memory of Joseph
Donin, Haifa, July 2004.

\section{Left and right skew-differential operators}
\label{sec-lrdif}

 Let $k$ be a field of characteristic zero, $G$ an abelian group,
and $V\in {}^{kG}_{kG}\YDcat $ a finite dimensional \YD module
with completely reducible $kG$-action. Let $\lcoa :V\to kG\ot V$
and $\lact :kG \ot V\to V$ denote the left coaction and left
action of $kG$ on $V$, respectively. Set $n:=\dim _kV$. By
assumption there exist nonzero numbers $q_{ij}\in k$, a basis
$\{x_i\,|\,1\le i\le n\}$ of $V$, and for each $i$ with $1\le i\le
n$ an element $g_i\in G$ such that
\begin{align}\label{eq-YDG}
g_i\lact x_j=&q_{ij}x_j,& \lcoa (x_j)=&g_j\ot x_j
\end{align}
for all $i,j\in \{1,2,\ldots ,n\}$. Then the braiding $\sigma \in
\End _k(V\ot V)$ of $V$ where
\begin{align*}
\sigma (v\ot w)=&(v_{(-1)}\lact w)\ot v_{(0)},& \sigma ^{-1}(v\ot
w)=&w_{(0)}\ot (\antip ^{-1}(w_{(-1)})\lact v),
\end{align*}
and $\lcoa (v)=v_{(-1)}\ot v_{(0)}$ for $v\in V$, is \textit{of
diagonal type}. Let $\cB (V)$ denote the Nichols algebra generated
by $V$. More precisely, as proved in \cite{a-Schauen96} and noted
in \cite[Prop.~2.11]{inp-AndrSchn02},
\begin{align*}
\cB (V)=k\oplus V\oplus \bigoplus _{m=2}^\infty V^{\ot m}/\ker S_m
\end{align*}
where $S_m\in \End _k(V^{\ot m})$, $S_{1,j}\in \End _k(V^{\ot
j+1})$,
\begin{align*}
S_m&=\prod _{j=1}^{m-1}(\id ^{\ot m-j-1}\ot S_{1,j})
=\prod _{j=1}^{m-1}(S_{j,1}\ot \id ^{\ot m-j-1}),\\
S_{1,j}&=\id +\sigma ^{-1}_{12}+\sigma ^{-1}_{12}\sigma ^{-1}_{23}
+\cdots +\sigma ^{-1}_{12}\sigma ^{-1}_{23}\cdots \sigma
^{-1}_{j,j+1},\\
S_{j,1}&=\id +\sigma ^{-1}_{j,j+1}+\sigma ^{-1}_{j,j+1}\sigma
^{-1}_{j-1,j} +\cdots +\sigma ^{-1}_{j,j+1}\cdots \sigma
^{-1}_{2,3}\sigma ^{-1}_{12}
\end{align*}
(in leg notation) for $m\ge 2$ and $j\in \ndN _0$. The maps $S_m$
are called \textit{braided symmetrizer}. They have analogues where
the braiding $\sigma $ is used instead of $\sigma ^{-1}$. Set $\cB
(V)_i:=V^{\ot i}/\ker S_i\subset \cB (V)$. The algebra $\cB (V)$
is $\ndZ ^n$-graded, where the degrees of the generators $x_i$ are
$\deg \,x_i=\Ndbasis _i$, and $\{\Ndbasis _i\,|\,1\le i\le n\}$ is
a basis of the $\ndZ $-module $\ndZ ^n$.

Note that in our setting $V$ is additionally a \YD module over
$k\ndZ ^n$ where
\begin{align}\label{eq-YDZn}
e_i\lactZ x_j=&q_{ji}x_j,& \lcoaZ (x_j)=e_j\ot x_j.
\end{align}
Here $\{e_i\,|\,1\le i\le n\}$ is a fixed basis of the $\ndZ
$-module $\ndZ ^n$, and $\lactZ $ and $\lcoaZ $ denote the left
action and left coaction of $k\ndZ ^n$ on $V$. In order to avoid
misunderstandings we will use the exponential notation for
elements of $\ndZ ^n$, thus $e_i^{-1}$ is the inverse of $e_i$ in
the group $\ndZ ^n$. The braiding $\sigma $ of $V$ commutes both
with the action of $g_i$ and the action of $e_i$: for all
$i,j,m\in \{1,2,\ldots ,n\}$ one has
\begin{align}\label{eq-gesigmacomm}
\sigma (g_i\lact (x_j\ot x_m))=&g_i\lact \sigma (x_j\ot x_m),&
\sigma (e_i\lactZ (x_j\ot x_m))=&e_i\lactZ \sigma (x_j\ot x_m).
\end{align}

The dual vector space $V^*$ gives rise to \textit{left} and
\textit{right skew-differential operators} $\ldif _i$ and $\rdif
_i$ on $\cB (V)$, where $1\le i\le n$, in the following way. Let
$\{y_i\,|\,1\le i\le n\}$ denote the basis of $V^*$ dual to
$\{x_i\,|\,1\le i\le n\}$. For $m\in \ndN $, $\rho \in \cB (V)_m$
and $i\in \{1,2,\ldots ,n\}$ set
\begin{align}\label{eq-lrdif}
\ldif _i(1)=&\rdif _i(1)=0,& \ldif _i(\rho )=&\rho '_i,&
\rdif_i(\rho )=&\rho ''_i,
\end{align}
where
\begin{align*}
S_{1,m-1}(\rho )=&\sum _{l=1}^n x_l\ot \rho '_l,&
 S_{m-1,1}(\rho )=&\sum _{l=1}^n \rho ''_l\ot x_l.
\end{align*}
Note that the antipode $\kappa $ of group algebras satisfies
$\antip ^2=\id $ and hence left and right duals of \YD modules
coincide. One can consider the vector space $\Lin _k\{\ldif
_i\,|\,1\le i\le n\}$ as a \YD module over $kG$ dual to $V$, that
is for $i,j\in \{1,2,\ldots ,n\}$ one has
\begin{align}
 g_i\lact \ldif _j=&q_{ij}^{-1}\ldif _j,&
 \lcoa (\ldif _j)=g_j^{-1}\ot \ldif _j.
\end{align}
Similarly, the vector space $\Lin _k\{\rdif _i\,|\,1\le i\le n\}$
becomes a \YD module over $k\ndZ ^n$ dual to $V$, that is one has
\begin{align}
 e_i\lactZ \rdif _j=&q_{ji}^{-1}\rdif _j,&
 \lcoaZ (\rdif _j)=e_j^{-1}\ot \rdif _j
\end{align}
for all $i,j\in \{1,2,\ldots ,n\}$. Moreover, the
skew-differential operators $\ldif _i$ and $\rdif _i$, where $1\le
i\le n$, and the braiding $\sigma $ satisfy the equations
\begin{align}
\begin{aligned}
 (\ldif _i\ot \id )(\sigma ^{-1}(x_j\ot x_m))=&(g_i^{-1}\lact x_j)\,
\ldif _i(x_m),\\
 (\id \ot \rdif _i)(\sigma ^{-1}(x_j\ot x_m))=&\rdif _i(x_j)\,e_i^{-1}
\lactZ x_m
\end{aligned}
\end{align}
for all $j,m\in \{1,2,\ldots ,n\}$. Therefore equations
(\ref{eq-lrdif}) give that
\begin{align}\label{eq-coprdif}
\begin{aligned}
 \ldif _i(\rho _1\rho _2)=&\ldif _i(\rho _1)\rho _2+g_i^{-1}\lact
\rho _1\,\ldif _i(\rho _2),\\
 \rdif _i(\rho _1\rho _2)=&\rho _1\,\rdif _i(\rho _2)
+\rdif _i(\rho _1)\,e_i^{-1}\lactZ \rho _2
\end{aligned}
\end{align}
for all $\rho _1,\rho _2\in \cB (V)$ and $i\in \{1,2,\ldots ,n\}$.
Note that by definition the skew-differential operators $\ldif _i$
and $\rdif _j$ commute:
\begin{align}\label{eq-difcomm}
 \ldif _i(\rdif _j(\rho ))=\rdif _j(\ldif _i(\rho )) \quad
 \text{for all $i,j\in \{1,2,\ldots ,n\}$, $\rho \in \cB (V)$}.
\end{align}
Moreover, the skew-differential operators $\ldif _i$ and $\rdif
_i$ and the group-like elements $g_j,e_j$, where $i,j\in
\{1,2,\ldots ,n\}$, satisfy the equations
\begin{align}\label{eq-gedifcomm}
\begin{aligned}
 g_j\lact (\ldif _i (g_j^{-1}\lact \rho ))=&q_{ji}^{-1}\ldif _i(\rho
),&
 e_j\lactZ (\ldif _i (e_j^{-1}\lactZ \rho ))=&q_{ij}^{-1}\ldif
_i(\rho ),\\
 g_j\lact (\rdif _i (g_j^{-1}\lact \rho ))=&q_{ji}^{-1}\rdif _i(\rho
),&
 e_j\lactZ (\rdif _i (e_j^{-1}\lactZ \rho ))=&q_{ij}^{-1}\rdif
_i(\rho )
\end{aligned}
\end{align}
for all $\rho \in \cB (V)$. It is a standard fact that the algebra
generated by the skew-differential operators $\ldif _i$, where
$i\in \{1,2,\ldots ,n\}$, is isomorphic to $\cB (V^*)$. Further,
as noted in \cite[Sect.~2.1]{a-Heck04a}, the assignment
$x_i\mapsto \ldif _i$ extends uniquely to an algebra isomorphism
$\iota :\cB (V)\to \cB (V^*)$. Finally, by Lemma 1 and Corollary 2
in \cite{a-Heck04a} there exists a unique bilinear map
$\paar{\cdot }{\cdot }:(\cB (V^*)\#kG)\ot \cB (V)\to \cB (V)$
which defines a $\cB (V^*)\#kG$-module algebra structure on $\cB
(V)$ extending the action of $\ldif _i$ and $g_i$ for all $i\in
\{1,2,\ldots ,n\}$.

\section{Finiteness conditions}
\label{sec-fini}

Retain the notation from the previous section. By a theorem of
Kharchenko \cite[Theorem~2]{a-Khar99} the algebra $\cB (V)$ has a
(restricted) Poincar\'e--Birkhoff--Witt basis consisting of
homogeneous elements with respect to the $\ndZ ^n$-grading of $\cB
(V)$. Let $\proots (\cB (V))\subset \ndZ ^n$ denote the set of
degrees of the (restricted) Poincar\'e--Birkhoff--Witt generators
of $\cB (V)$, counted with multiplicities. By the definition of
the $\ndZ ^n$-degree of $\cB (V)$ one clearly has $\proots (\cB
(V))\subset \ndN _0^n$. Set $\roots (\cB (V))=\proots (\cB
(V))\cup -\proots (\cB (V))$.

Let's consider the following finiteness conditions on $\cB (V)$.
\begin{itemize}
\item[(F1)] $\dim _k \cB (V)<\infty $,
\item[(F2)] the set $\proots (\cB (V))$ is finite,
\item[(F3)] $\mathrm{Dim}_k\,\cB (V)<\infty $,
\end{itemize}
where $\mathrm{Dim}_k$ denotes Gel'fand--Kirillov dimension.
Obviously one has the implications (F1)$\Rightarrow
$(F2)$\Rightarrow $(F3). Further, the condition (F1) holds if and
only if (F2) is satisfied and the heights of all restricted
Poincar\'e--Birkhoff--Witt generators of $\cB (V)$ are finite. It
is not known whether (F3) implies (F2).

As in \cite{inp-AndrSchn02} define
\begin{align}
\ad _\sigma x_i(\rho):=&x_i\rho-(g_i\lact \rho )x_i
\end{align}
for $\rho \in \cB (V)$ and $i\in \{1,2,\ldots ,n\}$. Consider the
sets
\begin{align*}
M_{i,j}:=\{(\ad _\sigma x_i)^m(x_j)\,|\,m\in \ndN _0\}
\end{align*}
for $i,j\in \{1,2,\ldots ,n\}$, $i\not=j$. By \cite[Lemma
20]{a-Rosso98}, if one assumes that (F3) holds then all $M_{i,j}$
are finite sets. More generally, by \cite[Lemma
3.7]{inp-AndrSchn02} or \cite[Sect.~4.1]{a-Heck04b} for given
$i,j\in \{1,2,\ldots ,n\}$ with $i\not=j$ the number
\begin{align}\label{eq-mij}
 m_{ij}:=\min\{m\in \ndN _0\,|\,\qnum{m+1}{q_{ii}}(q_{ii}^m
 q_{ij}q_{ji}-1)=0\}
\end{align}
is well-defined if and only if $M_{i,j}$ is finite. In this case
one has the relations $(\ad _\sigma x_i)^{m_{ij}+1}(x_j)=0$ and
$(\ad _\sigma x_i)^{m_{ij}}(x_j)\not=0$.

Fix $i\in \{1,2,\ldots ,n\}$. By the above paragraph, if $M_{i,j}$
is finite for all $j\in \{1,2,\ldots ,n\}$ with $j\not=i$ (for
example if (F3) holds) then one can introduce a $\ndZ $-linear
mapping $s_i:\ndZ ^n\to \ndZ ^n$ as follows:
\begin{align}\label{eq-si}
s_i(\Ndbasis _j):=&
 \begin{cases}
 -\Ndbasis _i & \text{if $j=i$,}\\
 \Ndbasis _j+m_{ij}\Ndbasis _i & \text{if $j\not=i$.}
 \end{cases}
\end{align}
In \cite[Ch.\,5, \S2]{b-BourLie4-6} such maps are called
pseudo-reflections. Note that $s_i^2=\id $.

\section{Transformations of Nichols algebras}
\label{sec-trafo}

Assume for the whole section that $\cB (V)$ is a rank $n$ Nichols
algebra of diagonal type. Further, suppose that $i\in \{1,2,\ldots
,n\}$ is chosen such that $M_{i,j}$ is finite for all $j\in
\{1,2,\ldots ,n\}$ with $j\not=i$. We describe the construction of
a Nichols algebra associated to $i$.

Since $k\rdif _i\in {}_{k\ndZ }^{k\ndZ }\YDcat $ one can construct
the smash product
\begin{align*}
H_i:=k[\rdif _i] \#k[e_i,e_i^{-1}].
\end{align*}
It has a unique Hopf algebra structure satisfying the formulas
\begin{align}\label{eq-HiHopf}
e_i\rdif _i=&q_{ii}^{-1}\rdif _ie_i,& \copr (e_i)=&e_i\ot e_i,&
\copr (\rdif _i)=1\ot \rdif _i+\rdif _i\ot e_i^{-1}.
\end{align}
By Equations (\ref{eq-gedifcomm}) and (\ref{eq-coprdif}) one
obtains that $\cB (V)$ is an $H_i$-module algebra, where $e_i$ and
$e_i^{-1}$ act via $\lactZ $ and $\rdif _i$ acts by evaluation. By
slight abuse of notation the symbol $\lactZ $ will also be used
for the left action of $H_i$ on $\cB (V)$. Let $(\cB (V)\op
\#H\cop _i)\op $ denote the opposite algebra of the smash product
of $\cB (V)\op $ and $H\cop _i$. Recall that it contains $\cB (V)$
and $H\opcop _i$ as subalgebras and one has
\begin{align}\label{eq-opsmash}
\rho h=&h_{(1)}(h_{(2)}\lactZ \rho ),
\end{align}
and in particular
\begin{align}\label{eq-opsmashyi}
\rho \rdif _i=&\rdif _i\cdot (e_i^{-1}\lactZ \rho )+\rdif _i(\rho
),
\end{align}
for all $\rho \in \cB (V)$ and $h\in H_i$, where $\copr
(h)=h_{(1)}\ot h_{(2)}$ denotes the coproduct of $h\in H_i$.
Further, $(\cB (V)\op \#H\cop _i)\op $ is a \YD module over $kG$
where the left action $\lacti $ and left coaction $\lcoai $ on
$(\cB (V)\op \#H\cop _i)\op $ are given by
\begin{align}\label{eq-YDBext}
\begin{aligned}
 \lcoai (x_j)=&g_j\ot x_j,&
 \lcoai (\rdif _i)=&g_i^{-1}\ot \rdif _i,&
 \lcoai (e_i)=&1\ot e_i,\\
 g_j\lacti x_m=&q_{jm}x_m,&
 g_j\lacti \rdif _i=&q_{ji}^{-1}\rdif _i,&
 g_j\lacti e_i=&e_i
\end{aligned}
\end{align}
for all $j,m\in \{1,2,\ldots ,n\}$.

\begin{satz}\label{s-trafo}
Assume that $i\in \{1,2,\ldots ,n\}$ such that $M_{i,j}$ is finite
for all $j\in \{1,2,\ldots ,n\}$ with $j\not=i$. Let $V_i$ denote
the subspace of the algebra $(\cB (V)\op \#H\cop _i)\op $
generated by the set $\{(\ad _{\sigma
}x_i)^{m_{ij}}(x_j)\,|\,j\not=i\}\cup \{\rdif _i\}$.  The
subalgebra $\cB _i$ of $(\cB (V)\op \#H\cop _i)\op $ generated by
$V_i$ is isomorphic to the Nichols algebra $\cB (V_i)$, and the
relation $\proots (\cB _i)=(s_i(\proots (\cB (V)))\setminus
\{-\Ndbasis _i\})\cup \{\Ndbasis _i\}$ holds.
\end{satz}

\begin{bem} It will be clear from the construction that the
transformation in Proposition \ref{s-trafo} is invertible. This
fact is in accord with the relation $s_i^2=\id $.
\end{bem}

\begin{bew}
One easily checks that $V_i$ is a \YD submodule of $(\cB (V)\op
\#H\cop _i)\op $ over $kG$ where the \YD module structure of $(\cB
(V)\op \#H\cop _i)\op $ is described above the proposition.

In the following we will sometimes consider $(\cB (V)\op \#H\cop
_i)\op $ as a vector space and as such we identify it with $\cB
(V)\ot H_i$.

We show that
\begin{align}\label{eq-Bi}
 \cB _i=\ker \ldif _i\ot k[\rdif _i](\subset \cB (V)\ot H_i).
\end{align}
Recall that $\rdif _i$ may have finite order, but it is always the
same as the order of $x_i$. Note also that by Equations
(\ref{eq-HiHopf}), (\ref{eq-opsmash}), and (\ref{eq-difcomm}) the
space $\ker \ldif _i\ot k[\rdif _i]$ is a subalgebra of $(\cB
(V)\op \#H\cop _i)\op $. The inclusion $\cB _i\subset \ker \ldif
_i\ot k[\rdif _i]$ holds by definition of $\cB _i$. In order to
prove that $\cB _i\supset \ker \ldif _i\ot k[\rdif _i]$ one first
checks that $\cB (V)\cong \ker \ldif _i\ot k[x_i]$ (as graded
vector spaces) using methods similar to \cite[Lemma 2.2]{a-Jos77}.
Further, $\ker \ldif _i$ is generated as an algebra by $\bigcup
_{j\,|\,j\not=i}M_{i,j}$. Now one proves by induction that
\begin{align}
 \rdif _i((\ad _\sigma x_i)^m(x_j))=&q_{ji}^{-1}\qnum{m}{q_{ii}^{-1}}
 (1-q_{ii}^{m-1}q_{ij}q_{ji})(\ad _\sigma x_i)^{m-1}(x_j)
\end{align}
for all $m\in \ndN _0$ and all $j\in \{1,2,\ldots ,n\}$. Hence by
Equation (\ref{eq-opsmashyi}) one obtains that $\ker \ldif _i\ot
k[\rdif _i]$ is generated by $V_i$.

In order to prove that the subalgebra $\cB _i$ of $(\cB (V)\op
\#H\cop _i)\op $ is a Nichols algebra it suffices to find $n$
appropriate skew-differential operators on $\cB _i$. More
precisely, according to the fact that the \YD module $V_i$
generates $\cB _i$, one has to find for all $m\in \{1,2,\ldots
,n\}$ a map $Y_m\in \End (\cB _i)$ such that
\begin{align}\label{eq-Ycond}
\begin{gathered}
\begin{aligned}
 Y_j(\rho _1\rho _2)=&Y_j(\rho _1)\rho _2
 +(g_i^{-m_{ij}}g_j^{-1}\lact \rho _1)\,Y_j(\rho _2),\\
 Y_i(\rho _1\rho _2)=&Y_i(\rho _1)\rho _2
 +(g_i\lact \rho _1)\,Y_i(\rho _2),
\end{aligned}\\
 Y_l((\ad _\sigma x_i)^{m_{ij}}(x_j))=\delta _{lj},\quad
 Y_l(\rdif _i)=\delta _{li} \quad \text{(Kronecker's delta)}
\end{gathered}
\end{align}
for all $\rho _1,\rho _2\in \cB _i$, and $j,l\in \{1,2,\ldots
,n\}$, $j\neq i$. For $\rho \in \cB _i$ set
\begin{align*}
 Y_i(\rho ):=&\ad _\sigma x_i(\rho)=x_i\rho -(g_i\lact \rho )x_i.
\end{align*}
By definition of $\ad _\sigma $, $Y_i$ satisfies the second
equation of (\ref{eq-Ycond}) for all $\rho _1,\rho _2 \in (\cB
(V)\op \#H\cop _i)\op $. Moreover, by Equation (\ref{eq-coprdif})
one obtains for all $\rho \in \ker \ldif _i$ the formula
\begin{align*}
 \ldif _i(x_i\rho -(g_i\lacti \rho )x_i)=\rho -g_i^{-1}\lacti
 (g_i\lacti \rho )=0.
\end{align*}
Hence $Y_i$ maps $\ker \ldif _i\subset \cB (V)$ onto itself. By
Equations (\ref{eq-YDBext}) and (\ref{eq-opsmashyi}) one gets also
\begin{align*}
 Y_i(\rdif _i)=&x_i\rdif _i-(g_i\lacti \rdif _i)x_i
 =x_i\rdif _i-q_{ii}^{-1}\rdif _ix_i=\rdif _i(x_i)=1.
\end{align*}
Finally, the equation $Y_i((\ad _\sigma x_i)^{m_{ij}}(x_j))=(\ad
_\sigma x_i)^{m_{ij}+1}(x_j)=0$ holds for all $j\in \{1,2,\ldots
,n\}$, $j\not=i$, by the definition of $m_{ij}$.

Let's now turn to the construction of $Y_j$, where $j\in
\{1,2,\ldots ,n\}$, $j\not=i$. Set $\lambda _j:=\paar{\iota ((\ad
_\sigma x_i)^{m_{ij}}(x_j))}{(\ad _\sigma x_i)^{m_{ij}}(x_j)}$
(for the notation see the last paragraph of Section
\ref{sec-lrdif}). By \cite[Sect.~4.1]{a-Heck04b} and the
definition of $m_{ij}$ one obtains that $\lambda _j\in k\setminus
\{0\}$. Define
\begin{align*}
 Y_j(\rho \ot (\rdif _i)^m):= \frac{1}{\lambda _j}
 \paar{\iota ((\ad _\sigma x_i)^{m_{ij}}(x_j))}{\rho }
 \ot (\rdif _i)^m
\end{align*}
for all $\rho \in \ker \ldif _i$, $m\in \ndN _0$, and $j\in
\{1,2,\ldots ,n\}$, $j\not=i$. By \cite[Eqn.~(2)]{a-Heck04b} the
coproduct of $\iota ((\ad _\sigma x_i)^{m_{ij}}(x_j))\in \cB
(V^*)\#kG $ takes the form
\begin{align*}
\iota ((\ad _\sigma x_i)^{m_{ij}}(x_j))\ot 1+ \sum _{m=0}^{m_{ij}}
c_i(\ldif _i)^{m_{ij}-m}g_i^{-m}g_j^{-1}\ot \iota ((\ad _\sigma
 x_i)^{m}(x_j))
\end{align*}
for certain $c_i\in k$, where $c_{m_{ij}}=1$. Since $\ldif _i$
vanishes on $\ker \ldif _i$, the first equation of
(\ref{eq-Ycond}) holds for all $\rho _1\in \ker \ldif _i$, $\rho
_2\in \cB _i$. It remains to show that the first equation of
(\ref{eq-Ycond}) is valid for $\rho _1=\rdif _i$ and $\rho _2\in
\ker \ldif _i$. By Equations (\ref{eq-opsmashyi}) and
(\ref{eq-difcomm}) one obtains that
\begin{align}
Y_j(\rdif _i\,\rho_2)=&Y_j((e_i\lactZ \rho _2)\ot \rdif _i - \rdif
_i(e_i\lactZ \rho _2))\notag \\
\label{eq-ni1}
 =&Y_j(e_i\lactZ \rho _2)\ot \rdif _i - \rdif _i(Y_j(e_i\lactZ \rho
 _2)).
\end{align}
On the other hand, by Equations (\ref{eq-opsmashyi}) and
(\ref{eq-gedifcomm}) one gets
\begin{align*}
(g_i^{-m_{ij}}g_j^{-1}\lact \rdif _i)\,Y_j(\rho _2)=&
q_{ii}^{m_{ij}}q_{ji}\Big( (e_i\lactZ Y_j(\rho _2))\ot \rdif _i
-\rdif _i(e_i\lactZ Y_j(\rho _2))\Big).
\end{align*}
By Equation (\ref{eq-gedifcomm}) the latter coincides with
(\ref{eq-ni1}).

Finally, it has to be shown that $\proots (\cB _i)=(s_i(\proots
(\cB (V)))\setminus \{-\Ndbasis _i\})\cup \{\Ndbasis _i\}$. In the
$\ndZ ^n$-graded algebra $(\cB (V)\op \#H\cop _i)\op $ the
elements $x_j$ and $\rdif _i$ have degree $\Ndbasis _j$ and
$-\Ndbasis _i$, respectively, for all $j\in \{1,2,\ldots ,n\}$.
Hence the elements $(\ad _\sigma x_i)^{m_{ij}}(x_j)$ have degree
$\Ndbasis _j+m_{ij}\Ndbasis _i=s_i(\Ndbasis _j)$ for all $j\in
\{1,2,\ldots ,n\}$, $j\not=i$. Fix the $\ndZ ^n$-degrees of the
generators of $\cB _i$ by
\begin{align}
\deg \,\rdif _i:=&\Ndbasis _i,& \deg \,(\ad _\sigma
x_i)^{m_{ij}}(x_j):=\Ndbasis _j
\end{align}
for $j\in \{1,2,\ldots ,n\}$, $j\not=i$. Then $s_i(\proots
(\cB _i))$ is exactly the set of degrees of the (restricted)
Poincar\'e--Birkhoff--Witt generators of $\cB _i$ in $(\cB (V)\op
\#H\cop _i)\op $.

Recall that any Nichols algebra $\cB $ of rank $n$ is isomorphic
as a $\ndZ ^n$-graded vector space to the algebra $k[x_r\,|\,r\in
\proots (\cB )]/(x_r^{h_r}\,|\,r\in \proots (\cB ),h_r<\infty )$.
Here $h_r$ denotes the height of the (restricted)
Poincar\'e--Birkhoff-Witt generator corresponding to $r\in \proots
(\cB )$, and by Kharchenkos theorem it is uniquely determined by
the $\ndZ ^n$-degree of $x_r$. Therefore it suffices to know the
multiplicities of the $\ndZ ^n$-homogeneous components of $\cB $
in order to determine $\proots (\cB (V))$. This fact and the
equation $\cB (V)\cong \ker \ldif _i\ot k[x_i]$ allows us to
conclude that (with obvious interpretation) $\proots (\ker \ldif
_i\ot k[\rdif _i])=(\proots (\cB (V))\setminus \{\Ndbasis
_i\})\cup \{-\Ndbasis _i\}$. Since the determination of $\proots $
can be performed using different total orderings of $\ndZ ^n$, we
conclude from Equation (\ref{eq-Bi}) that $s_i(\proots (\cB _i))=
(\proots (\cB (V))\setminus \{\Ndbasis _i\})\cup \{-\Ndbasis
_i\}$. This proves the proposition.
\end{bew}

\begin{bem}
1. Proposition \ref{s-trafo} has the following interpretation. Set
$\roots :=\roots (\cB (V))\subset \ndZ ^n$. Then the set $\roots
(\cB _i)$ coincides with $\roots $ with respect to the basis
$\{s_i(\Ndbasis _j)\,|\,1\le j\le n\}$ of $\ndZ ^n$. With other
words, the transformation doesn't change $\roots (\cB (V))$, it
changes only the basis of $\ndZ ^n$. Additionally, this base
change is performed in such a way that the new basis is a subset
of $\roots (\cB (V))$.

2. The constants $m_{ij}$ appearing in the definition of the map
$s_i$, where $j\in \{1,2,\ldots ,n\}$ and $j\not=i$, depend on the
structure constants of the braiding. The latter usually change if
one performs a transformation. However in a special case, namely
if the braiding is of Cartan type, this change is not essential.
This situation is the one which is best understood, and it will
also be analyzed in more detail in Section \ref{sec-Cartan}.
\end{bem}

\section{The Weyl--Brandt groupoid associated to a Nichols
algebra} \label{sec-WBg}

Let $\cG $ be a nonempty set, $D\subset \cG \times \cG $ a
nonempty subset, and $\circ :D\to \cG $ a map of sets. The pair
$(\cG ,\circ )$ is called \textit{Brandt groupoid} if it satisfies
the following conditions (see for example \cite[Sect.
3.3]{b-ClifPres61}).
\begin{itemize}
\item If $(x,y)\in D$ then each of the three elements $x,y,x\circ
y$ is uniquely determined by the other two.
\item If $(x,y),(y,z)\in D$ then $(x\circ y,z),(x,y\circ z)\in D$
and $(x\circ y)\circ z=x\circ (y\circ z)$.
\item If $(x,y),(x\circ y,z)\in D$ then $(y,z),(x,y\circ z)\in D$
and $(x\circ y)\circ z=x\circ (y\circ z)$.
\item If $(y,z),(x,y\circ z)\in D$ then $(x,y),(x\circ y,z)\in D$
and $(x\circ y)\circ z=x\circ (y\circ z)$.
\item For all $x\in \cG $ there exist unique elements $e,f,y\in \cG$
such that $(e,x)$, $(x,f),(y,x)\in D$, $e\circ x=x\circ f=x$, and
$y\circ x=f$.
\item If $e\circ e=e$, $f\circ f=f$ for certain $e,f\in \cG $ then
there exists $x\in \cG $ such that $e\circ x=x\circ f=x$.
\end{itemize}

Let $\cB (V)$ be a rank $n$ Nichols algebra of diagonal type. Let
$E_0:=(\Ndbasis _1,\ldots ,\Ndbasis _n)$ denote the (ordered)
standard basis of $\ndZ ^n$. Define
\begin{align*}
W(V):=\{(s,E)\,|\,&\text{$s\in \Aut (\ndZ ^n)$, $E$ is an ordered
basis of $\ndZ ^n$,}\\
&\text{there exist $m_1,m_2\in \ndN _0$, $m_1\le m_2$},\\
&\text{and $i_1,\ldots ,i_{m_2}\in \{1,2,\ldots ,n\}$, such
that}\\
&\text{$s_{i_m}\ldots s_{i_1}(E_0)$ is well defined for all $m\le
m_2$ and }\\
&\text{$s_{i_{m_1}}\ldots s_{i_1}(E_0)=E$, $s=s_{i_{m_2}}\ldots
s_{i_{m_1+1}}$\}}.
\end{align*}
Note that $W(V)$ is not empty since $(\id ,E_0)\in W(V)$. Recall
that $s_i^2=\id $ whenever $s_i$ is defined. Hence there is a
natural Brandt groupoid structure on $\cG=W(V)$ such that
$(s,E)\circ (t,F)$ is defined (and is then equal to $(st,F))$ if
and only if $t(F)=E$. We call $W(V)$ the \textit{Weyl--Brandt
groupoid of} $\cB (V)$.

The definition of $W(V)$ gives an important consequence of
Proposition \ref{s-trafo}.

\begin{folg}\label{f-finiteW}
If $\cB (V)$ is a rank $n$ Nichols algebra of diagonal type
satisfying (F2) then $W(V)$ is finite. In particular, the orbit of
any element and any ordered basis of $\ndZ ^n$ under the action of
$W(V)$ is finite.
\end{folg}

\section{Nichols algebras of Cartan type}
\label{sec-Cartan}

If the braiding $\sigma $ of $V$ is of diagonal type and the
structure constants $q_{ij}$ satisfy the equations
\begin{align}\label{eq-Cartantype}
q_{ij}q_{ji}=q_{ii}^{a_{ij}},\quad \text{$i,j\in \{1,2,\ldots
,n\}$,}
\end{align}
where $a_{ii}=2$ and $a_{ij}$ is a nonpositive integer for all
$i\not=j$, then one says that $\sigma $ (and $V$ and $\cB (V)$) is
of \textit{Cartan type}. It is then always assumed that the
$a_{ij}$ are maximal with the above properties. Note that if
$a_{ij}>0$ for some $j\not=i$ and $q_{ii}$ is not a root of unity
then by Rosso \cite[Lemma 20]{a-Rosso98} the Gel'fand--Kirillov
dimension of $\cB (V)$ is infinite.

\begin{lemma}\label{l-mij=-aij}
Assume that $V$ is an $n$-dimensional \YD module of Cartan type
and let $(a_{ij})_{i,j\in \{1,2,\ldots ,n\}}$ denote the
corresponding Cartan matrix.\\
(i) Suppose that $i\in \{1,2,\ldots ,n\}$ and $m\in \ndN _0$ such
that $m<-a_{ij}$ for at least one $j\in \{1,2,\ldots ,n\}$. Then
$q_{ii}^{m+1}\not=1$.\\
(ii) For the numbers in Equation (\ref{eq-mij}) one obtains
$m_{ij}=-a_{ij}$ for all $i,j\in \{1,2,\ldots ,n\}$.
\end{lemma}

\begin{bew}
To (i). Since $V$ is of Cartan type, one has
$q_{ii}^{a_{ij}}=q_{ij}q_{ji}$. Assume that $m+1\le -a_{ij}$ and
$q_{ii}^{m+1}=1$. One obtains that
$q_{ii}^{m+1+a_{ij}}=q_{ij}q_{ji}$, and $m+1+a_{ij}\le 0$. This is
a contradiction to the choice of $a_{ij}$.

To (ii). This follows from the definition of $m_{ij}$ and from
(i).
\end{bew}

One says that $\sigma $ is of \textit{finite type} if the Cartan
matrix $(a_{ij})_{i,j\in \{1,2,\ldots ,n\}}$ is of finite type.
There exist classification results of Rosso
\cite[Theorem~21]{a-Rosso98} and Andruskiewitsch and Schneider
\cite[Theorem~1.1]{a-AndrSchn00} on Nichols algebras of Cartan
type with finite Gel'fand--Kirillov dimension (F3) and finite
dimension (F1), respectively. The introduction of the Weyl--Brandt
groupoid $W(V)$ in the previous section allows to state a theorem
without technical assumptions on the numbers $q_{ij}$.

\begin{thm}\label{t-class}
Let $V$ be an $n$-dimensional \YD module of Cartan type with
corresponding Cartan matrix $C:=(a_{ij})_{i,j\in \{1,2,\ldots
,n\}}$.\\
(i) If $C$ is not of finite type then $\roots (\cB (V))$ is
infinite.\\
(ii) If $C$ is of finite type then $\roots (\cB (V))$ can be
identified with the set of roots of the semisimple Lie algebra
corresponding to $C$. Moreover, in any connected component the
heights of the (restricted) Poincar\'e--Birkhoff--Witt generators
depend only on the lengths of the roots corresponding to them.
\end{thm}

\begin{bew}
To (i). Assume that $\roots (\cB (V))$ is finite. By Proposition
\ref{s-trafo} for each $i\in \{1,2,\ldots ,n\}$ there exists a
Nichols algebra $\cB _i\cong \cB (V_i)$ of rank $n$. By Lemma
\ref{l-mij=-aij}(ii) one has $m_{ij}=-a_{ij}$. Choose the basis of
$V_i$ in such a way that the $j^\mathrm{th}$ basis vector is $(\ad
_\sigma x_i)^{-a_{ij}}(x_j)$ if $j\not=i$ and $\rdif _i$ if $j=i$,
respectively. Let $\{q(i)_{jm}\,|\,j,m\in \{1,2,\ldots ,n\}\}$
denote the set of structure constants of the braiding of $V_i$
with respect to this basis. By Equations (\ref{eq-YDBext}) one
gets
\begin{align*}
q(i)_{jm}=
\begin{cases}
q_{ii} & \text{if $j=m=i$,}\\
q_{ii}^{a_{im}}q_{im}^{-1}=q_{mi} & \text{if $j=i,m\not=i$,}\\
q_{ii}^{a_{ij}}q_{ji}^{-1}=q_{ij} & \text{if $j\not=i,m=i$,}\\
q_{ii}^{a_{ij}a_{im}}q_{im}^{-a_{ij}}q_{ji}^{-a_{im}}q_{jm}=
q_{ij}^{a_{im}}q_{im}^{-a_{ij}}q_{jm} & \text{if $j\not=i,m\not=i$.}\\
\end{cases}
\end{align*}
In particular, one obtains that
\begin{align}\label{eq-q(i)}
\begin{aligned}
 q(i)_{jj}=&q_{jj},& q(i)_{jm}q(i)_{mj}=&q(i)_{jj}^{a_{jm}}
\end{aligned}
\end{align}
for all $j,m\in \{1,2,\ldots ,n\}$. Hence $\cB _i$ is of the same
Cartan type as $\cB (V)$. This means that the linear maps $s$ in
the elements $(s,E)\in W(V)$ don't depend on the basis $E$ of
$\ndZ ^n$. By (F2) and Corollary \ref{f-finiteW} the Weyl--Brandt
groupoid $W(V)$ is finite. If $C$ is symmetrizable then $W(V)$ is
canonically isomorphic to $W\times B$ where $W$ is the Weyl group
associated to $C$ and $B$ is the orbit of $E_0$ under $W$. Here
$E_0:=(\Ndbasis _1,\ldots ,\Ndbasis _n)$ denotes the ordered
standard basis of $\ndZ ^n$. It is well-known
\cite[Ch.\,1,\,Theorem~4.8]{b-Hiller82} that $W$ is finite if and
only if the symmetrizable Cartan matrix is of finite type. Thus
for (i) it suffices to show that the group $W(V)$ is infinite
whenever $C$ is \textit{not} symmetrizable. Note that if $W(V')$
corresponding to a \YD submodule $V'$ of $V$ is infinite then
$W(V)$ is itself infinite. Further, if the Dynkin diagram
associated to a Cartan matrix is simply-laced or has no cycles
then it is symmetrizable. Thus we only have to show that $W(V)$ is
infinite if the corresponding Dynkin diagram is a cycle which is
not simply-laced. Further, if we remove a node from a cycle, then
we come again to the symmetrizable case. Thus it is sufficient to
consider cycles (which are not symmetrizable) such that after
removing an arbitrary node the resulting diagram is of finite
type. Using the classification result of Cartan matrices of finite
type one obtains easily (for similar argumentations confer also
\cite[Sect.~4.4]{a-AndrSchn00}) that such cycles have three nodes,
or the corresponding Cartan matrix is one of the following:
\begin{align*}
\begin{pmatrix}
 2&-2& 0&-1\\
-1& 2&-1& 0\\
 0&-1& 2&-2\\
-1& 0&-1& 2
\end{pmatrix},\qquad
\begin{pmatrix}
 2&-2& 0& 0&-1\\
-1& 2&-1& 0& 0\\
 0&-1& 2&-1& 0\\
 0& 0&-1& 2&-1\\
-1& 0& 0&-1& 2
\end{pmatrix}.
\end{align*}
In the last case Equations (\ref{eq-Cartantype}) imply that
$q_{11}^2=q_{22}=q_{33}=q_{44}=q_{55}=q_{11}$ and hence
$q_{12}q_{21}=q_{11}=1$. By the maximality assumption on $a_{ij}$
the latter means that this type of Cartan matrix does not appear.

Consider cycles with three nodes. The matrices $t_i$ of $s_i$,
$i\in \{1,2,3\}$, with respect to the basis $\{\Ndbasis
_1,\Ndbasis _2,\Ndbasis _3\}$ take the form
\begin{align*}
t_1=&\begin{pmatrix}
 -1 & \makebox[3ex]{$-a_{12}$} & \makebox[3ex]{$-a_{13}$}\\
  0 & 1 & 0\\
  0 & 0 & 1
\end{pmatrix},&
t_2=&\begin{pmatrix}
  1 & 0 & 0\\
 \makebox[3ex]{$-a_{21}$} & -1 & \makebox[3ex]{$-a_{23}$}\\
  0 & 0 & 1
\end{pmatrix},&
t_3=&\begin{pmatrix}
  1 & 0 & 0\\
  0 & 1 & 0\\
  \makebox[3ex]{$-a_{31}$} & \makebox[3ex]{$-a_{32}$} & -1
\end{pmatrix}.
\end{align*}
Without loss of generality one can suppose that $a_{12}<-1$. The
element $s_1s_2s_3$ (recall the independence of the $s_i$ from the
basis of $\ndZ ^n$) has the matrix
\begin{align*}
t_1t_2t_3{=}\begin{pmatrix}
a_{12}a_{21}{+}a_{13}a_{31}{-}1{-}a_{12}a_{23}a_{31} &
a_{12}{+}a_{13}a_{32}{-}a_{12}a_{23}a_{32} &
a_{13}{-}a_{12}a_{23}\\
a_{23}a_{31}-a_{21} & a_{23}a_{32}-1 & a_{23}\\
-a_{31} & -a_{32} & -1
\end{pmatrix}.
\end{align*}
Since we assumed (F2), by Corollary \ref{f-finiteW} this matrix
has to have finite order with respect to multiplication. This
means in particular that all eigenvalues of the (invertible)
matrix $t_1t_2t_3$ have to have absolute value 1, and hence its
trace is not bigger than 3. Further, if the trace is 3 then the
matrix has to be the identity. Since $a_{12}\le -2$, for the trace
of $t_1t_2t_3$ one obtains the relation
\begin{align*}
\mathrm{tr}(t_1t_2t_3)=a_{12}a_{21}+a_{13}a_{31}+a_{23}a_{32}
-a_{12}a_{23}a_{31}-3\ge 2{+}1{+}1{+}2{-}3=3.
\end{align*}
As $t_1t_2t_3$ is obviously not the identity this yields that
$s_1s_2s_3$ doesn't have finite order and hence $W(V)$ is
infinite. An analogous conclusion holds for the remaining cycle
with 4 nodes, where the matrix of $s_1s_2s_3s_4$ and its trace are
\begin{align*}
t_1t_2t_3t_4&=
\begin{pmatrix}
6 & 0 & 3 & -5\\ 3 & 0 & 1 & -2\\ 2 & 1 & 1 & -2\\ 1 & 0 & 1 & -1
\end{pmatrix},
& \mathrm{tr}(t_1t_2t_3t_4)=6>4.
\end{align*}

For the proof of (ii) we need the following lemma.

\begin{lemma}\label{l-DeltaC}
Let $C$ be a symmetrizable Cartan matrix of finite type with
corresponding Weyl group $W$. Let $\Delta _C$ denote the root
system corresponding to $C$. Let $\pi $ be a fixed set of simple
roots generating $\Delta _C$. Then
\begin{align*}
\{m\alpha \,|\,m\in \ndZ ,\alpha \in \Delta _C\}=\{\alpha \in \ndZ
\pi \, |\,W\alpha \subset \ndN _0\pi \cup -\ndN _0\pi \}.
\end{align*}
\end{lemma}

\begin{bew}
The inclusion $\subset$ in the equation of the lemma is
well-known. Further, if the rank $n$ of $C$ is one then the
inclusion $\supset $ is trivial.

Suppose that the inclusion $\supset $ in the lemma doesn't hold.
Without loss of generality the rank $n$ of $C$ is minimal with
this property. By the above remark one has $n\ge 2$. Let $\alpha
\in \ndZ \pi $ with $W\alpha \subset \ndN _0\pi \cup -\ndN _0\pi
$. Without loss of generality one can take $\alpha \in \ndN _0\pi
$. If $\alpha \notin \ndN \pi $ then $\alpha =m\beta $ for some
$\beta \in \Delta _C$, $m\in \ndZ $, by minimality of $n$.
Otherwise, since application of a simple reflection $s_i$ onto
$\alpha $ changes only 1 of its coefficients, $s_i(\alpha )\in
\ndN _0\pi $ by assumption on $\alpha $. Again, if $s_i(\alpha
)\notin \ndN \pi $ then $s_i(\alpha )=m\beta $ for some $\beta \in
\Delta _C$, $m\in \ndZ $, by minimality of $n$, and hence $\alpha
=ms_i(\beta )$. Finally, the case $W\alpha \subset \ndN \pi $ can
not appear since $w_0\alpha \in -\ndN \pi $ where $w_0$ is the
longest element of $W$.
\end{bew}

To (ii). If $V$ is of finite type then the Cartan matrix
$C=(a_{ij})_{i,j\in \{1,2,\ldots ,n\}}$ is symmetrizable and hence
$W(V)$ is isomorphic to $W\times B$. Define a $\ndZ $-linear map
$\phi :\ndZ \pi \to \ndZ ^n$, where $\pi =\{\alpha _1,\ldots
,\alpha _n\}$ is a fixed set of simple roots of the root system
$\Delta _C$ associated to $C$, by the formula
\begin{align*}
\phi (\alpha _i):=\Ndbasis _i.
\end{align*}
Note that $\phi $ commutes with the action of the maps $s_i$,
where $s_i$ are also interpreted as reflections on the set $\Delta
_C$ with respect to simple roots. Since all elements of $\Delta
_C$ may be obtained from simple roots by application of an element
of the Weyl group, one obtains $\roots (\cB (V))\supset \phi
(\Delta _C)$. Since the multiplicities of the degrees of the
generators of $\cB (V)$ are one, the multiplicity of $\phi
(\alpha)$ is one for all $\alpha \in \Delta _C$.

Assume now that $\alpha \in \roots (\cB (V))\setminus \phi (\Delta
_C)$. It is well-known that $\alpha \notin \ndZ \Ndbasis _i$ for
all $i\in \{1,2,\ldots ,n\}$. By application of elements in $W$
one obtains that $\alpha \notin \{m\phi (\beta )\,|\,m\in \ndZ
,\beta \in \Delta _C\}$. By Lemma \ref{l-DeltaC} there exists
$w\in W$ such that $w\alpha \notin \ndN _0\pi \cup -\ndN _0\pi $.
This is a contradiction to $\proots (\cB (V))\cup -\proots (\cB
(V))=\roots (\cB (V))=w(\roots (\cB (V)))$.

Recall that $q_{ii}=q_{jj}$ whenever $i,j\in \{1,2,\ldots ,n\}$
are in the same connected component, and any root $\alpha \in
\Delta _C$ can be written as $\alpha =w(\alpha _i)$ for some
simple root $\alpha _i$ and an element $w$ of the Weyl group. Thus
the last statement of the theorem follows from Equation
\ref{eq-q(i)}.
\end{bew}

\bibliographystyle{mybib}
\bibliography{quantum}

\end{document}